\pdfoutput=1
\RequirePackage{ifpdf}
\ifpdf 
\documentclass[pdftex]{sigma}
\else
\documentclass{sigma}
\fi

\usepackage[all]{xy}
\usepackage{lscape}

\numberwithin{equation}{section}

\newtheorem{Theorem}{Theorem}[section]
\newtheorem{Corollary}[Theorem]{Corollary}
\newtheorem{Proposition}[Theorem]{Proposition}
\newtheorem{Question}[Theorem]{Question}
 { \theoremstyle{definition}
\newtheorem{Definition}[Theorem]{Definition}
\newtheorem{Example}[Theorem]{Example}
\newtheorem{Remark}[Theorem]{Remark} }

\begin{document}

\newcommand{\arXivNumber}{1411.7063}

\allowdisplaybreaks

\renewcommand{\thefootnote}{$\star$}

\renewcommand{\PaperNumber}{062}

\FirstPageHeading

\ShortArticleName{Topological Monodromy of an Integrable Heisenberg Spin Chain}

\ArticleName{Topological Monodromy of an Integrable Heisenberg\\ Spin Chain\footnote{This paper is a~contribution to the Special Issue on
Algebraic Methods in Dynamical Systems.
The full collection is available at
\href{http://www.emis.de/journals/SIGMA/AMDS2014.html}{http://www.emis.de/journals/SIGMA/AMDS2014.html}}}

\Author{Jeremy LANE}

\AuthorNameForHeading{J.~Lane}

\Address{Department of Mathematics, University of Toronto,\\ 40 St.~George Street, Toronto, Ontario, Canada M5S 2E4}
\Email{\href{mailto:jeremy.lane@mail.utoronto.ca}{jeremy.lane@mail.utoronto.ca}}
\URLaddress{\url{https://www.math.toronto.edu/cms/lane-jeremy/}}

\ArticleDates{Received November 27, 2014, in f\/inal form July 29, 2015; Published online July 31, 2015}

\Abstract{We investigate topological properties of a completely integrable system on $S^2\times S^2 \times S^2$ which was recently shown to have a  Lagrangian f\/iber dif\/feomorphic to $\mathbb{R} P^3$ not displaceable by a Hamiltonian isotopy~[Oakley~J., Ph.D.~Thesis,
  University of Georgia, 2014]. This system can be viewed as integrating the determinant, or alternatively, as integrating a classical Heisenberg spin chain. We show that the system has non-trivial topological monodromy and relate this to the geometric interpretation of its integrals.}


\Classification{37J35; 53D12}

\renewcommand{\thefootnote}{\arabic{footnote}}
\setcounter{footnote}{0}

\section{Introduction}

The Hamiltonian
\begin{gather*}
  H(X,Y,Z) = \sqrt{3 + \langle X,Y\rangle + \langle Y,Z\rangle + \langle Z,X\rangle}
\end{gather*}
models pairwise interaction of three identical spin vectors $X,Y,Z\in S^2$ which are f\/ixed to the vertices of an equilateral triangle. Systems of this type -- called Heisenberg spin chains -- are of interest to physicists as they provide a classical model for quantum spin in a f\/ixed lattice.  Together with the Hamiltonians
\begin{gather*}
  I(X,Y,Z) = \langle X+Y+Z,e_3\rangle\qquad \mbox{and} \qquad J(X,Y,Z) = \det(X,Y,Z),
\end{gather*}
this Heisenberg spin chain becomes a completely integrable system on $(S^2\times S^2\times S^2, \omega_{\rm STD}\oplus \omega_{\rm STD} \oplus \omega_{\rm STD})$.  Lagrangian f\/ibers of this system were recently studied in \cite{oakley}, and it was shown that the system has a Lagrangian f\/iber $L\cong \mathbb{R} P^3$ which is not displaceable by Hamiltonian dif\/feomorphisms.

There is an analogous coupled spin system given by the Hamiltonians
\begin{gather*}
  H_1(X,Y) = \sqrt{1-\langle X,Y\rangle}\qquad \mbox{and} \qquad H_2(X,Y) = \langle X + Y,e_3\rangle,
\end{gather*}
which has provided interesting examples of non-displaceable Lagrangian tori in $(S^2 \times S^2, \omega_{\rm STD}\oplus\omega_{\rm STD})$ \cite{enpo, fooo-exotic-torus,oakley-usher,weiwei-wu}. In fact, $(H_1,H_2)$ is the moment map of a Hamiltonian $T^2$-action on the complement of the Lagrangian sphere $\widetilde\Delta$ of anti-diagonal elements~$(X,-X)$, where~$H_1$ fails to be smooth.

In this note we study the integrable system $\mathcal{H} = (H,I,J)$ from the perspective of Lagrangian torus f\/ibrations and global action-angle coordinates, as introduced by Duistermaat in~\cite{duistermaat}.  Unlike the system on $S^2\times S^2\setminus \widetilde\Delta$, the system $\mathcal{H}$ is not toric on the set where it is smooth, $S^2\times S^2 \times S^2 \setminus L$, and cannot be made so: there is a global obstruction to the existence of a dif\/feomorphism~$f$ such that the Hamiltonian f\/low of~$f\circ J$ is periodic. We compute this obstruction, called topological monodromy, explicitly
\begin{Theorem}\label{monodromy}
  The topological monodromy of the system $\mathcal{H}$ is generated by the matrix
  \begin{gather*}
      \left(\begin{matrix}
      1&0&0\\
      0&1&3\\
      0&0&1
     \end{matrix}\right).
  \end{gather*}
  Thus the system does not admit global action coordinates $($since $H$ and~$I$ are already periodic, this implies that there does not exist a diffeomorphism $f\colon \operatorname{Im}(J) \rightarrow \mathbb{R}$ such that the Hamiltonian $f\circ J$ is periodic$)$.
\end{Theorem}

\begin{figure}[t]
        \centering
\includegraphics{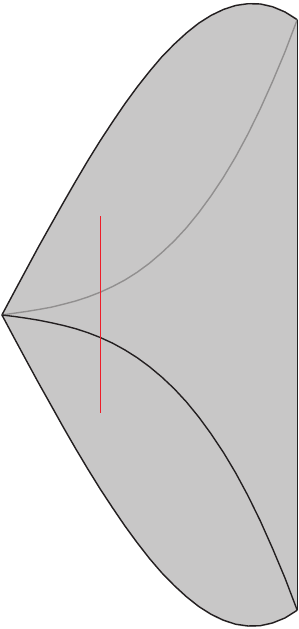}

                    \caption{The moment map image is a solid with one `orbifold' corner, four `toric' faces, and a~`critical line' through the interior (red).}
                    \label{figure 1}
\end{figure}

\begin{figure}[t]
        \centering
\includegraphics{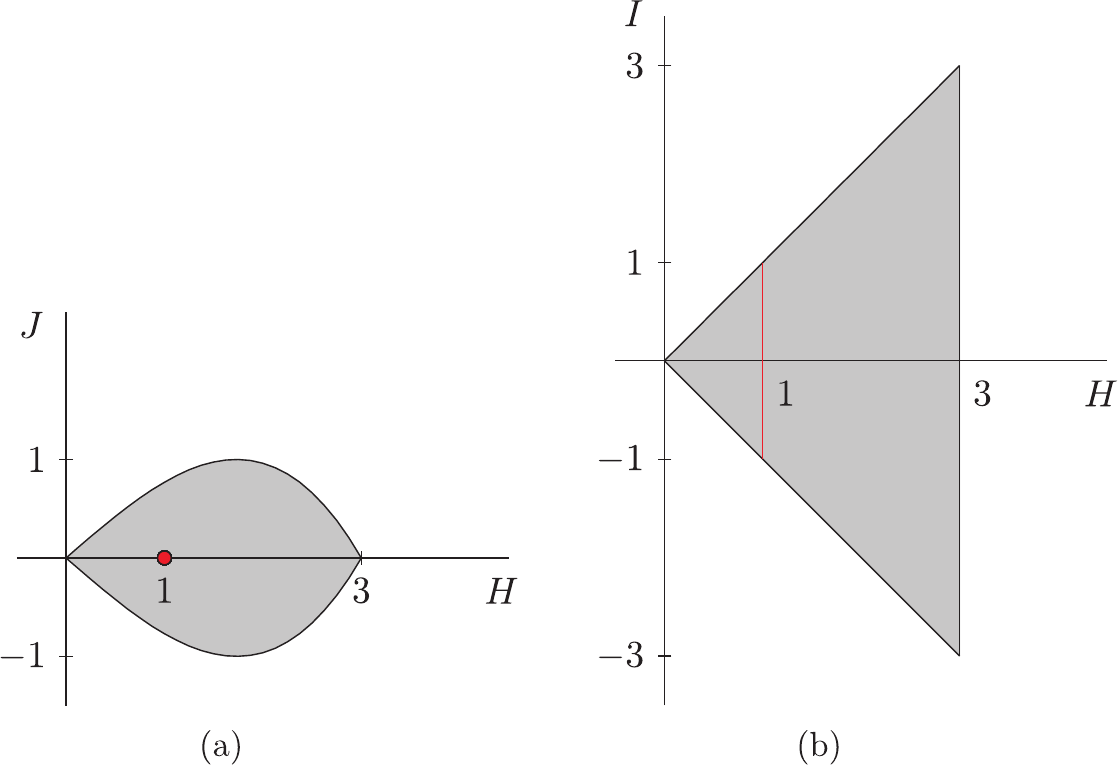}
        \caption{Projections of the moment map image.}
        \label{figure 2}
\end{figure}

Although topological monodromy of a given integrable system may be computed directly using elliptic integrals~-- as has been done in \cite{bates, duistermaat} and many other places~-- we derive this result from general theorems of Zung and Izosimov about the structure of integrable systems near singularities~\cite{izosimov, zung1}.  Section~\ref{section 2} reviews the theory of singularities for integrable systems developed in~\cite{zung1} and the relation of this theory to the topological monodromy obstruction.  Section~\ref{section 3} establishes basic facts and notation for adjoint orbits in Lie algebras which we use throughout the paper.  Section~\ref{section 4} introduces the system~$\mathcal{H}$ and provides an algebraic proof that the Hamiltonians~$H$,~$I$,~$J$ Poisson commute.  In order to apply the theorems of Zung and Isozimov to this system, we must describe the topology of the critical f\/ibers of the map~$\mathcal{H}$. In Section~\ref{section 5} we use the underlying Euclidean geometry of the system and the structure of the Lie algebra~$\mathfrak{so}(3)$ to completely describe the critical set, critical values and image of~$\mathcal{H}$,
\begin{Theorem}\label{image of moment map}
  The image of the moment map $\mathcal{H}$ is the set of points $(r,s,t)\in \mathbb{R}^3$ that satisfy the equations
  \begin{gather}\label{region}
     |t| \leq \sqrt{1 + 2\left(\frac{r^2-3}{6}\right)^3 - 3\left(\frac{r^2-3}{6}\right)^2}, \qquad |s| \leq r, \qquad \mbox{and}\qquad 0 \leq  r \leq 3
  \end{gather}
  $($see Figs.~{\rm \ref{figure 1}} and~{\rm \ref{figure 2})}. The set of critical values consists of the boundary of this set, together with the line segment $\{(1,s,0)\colon {-}1<s<1\}$.
\end{Theorem}

We show in Section~\ref{section 7} that the critical f\/ibers above the line segment $\{(1,s,0)\colon {-}1<s<1\}$, which we call the `critical line', are all topologically stable, rank~1, focus-focus (see Section~\ref{section 2} for def\/initions). Theorem~\ref{monodromy} follows from this by the theorems of Zung and~Izosimov.

Although it is not necessary in order to prove Theorem~\ref{monodromy}, we give a hands-on proof that the regular f\/ibers of~$\mathcal{H}$ are all connected in Section~\ref{section 6}. Combining our description of the moment map image, and results of~\cite{duistermaat}, the Lagrangian torus f\/ibration over the set of regular values of $\mathcal{H}$ is determined up to f\/iber-preserving symplectomorphism by Theorem~\ref{monodromy} (see Remark~\ref{chern}).

\section{Singularities and topological monodromy\\ of integrable systems}\label{section 2}

In this section we review the theory of singularities for integrable systems developed in~\cite{zung1} and its relation to topological monodromy. We make several simplifying assumptions (such as f\/iber connectedness) and the reader should note that Theorem~\ref{zung classification} is not stated in full generality. Historical background and further details can be found in~\cite{bolsinov, zung1}.

\begin{Definition}
  An \emph{integrable system} is a triple $(M,\omega,\mathcal{H})$ where $(M,\omega)$ is a symplectic manifold of dimension $2n$, and $\mathcal{H} = (H_1,\ldots ,H_n)\colon M\rightarrow \mathbb{R}^n$ is a smooth map such that
\begin{enumerate}\itemsep=0pt
  \item[1)] the Poisson bracket $\{H_i,H_j\} = 0$ for all $i$, $j$, and
  \item[2)] there is an open, dense subset of $M$ where the functions $H_1,\ldots ,H_n$ are all functionally independent (i.e., $dH_1 \wedge \cdots \wedge dH_n \neq 0$).
\end{enumerate}
\end{Definition}

We denote the image of $\mathcal{H}$ by $B$, the set of regular values by $B_{\rm reg}$, and preimage $\mathcal{H}^{-1}(B_{\rm reg})$ by $M_{\rm reg}$.

\begin{Definition}
  An integrable system $(M,\omega,\mathcal{H})$ admits \emph{global action coordinates} if there is a~dif\/feomorphism $a = (a_1, \ldots ,a_n)\colon B \rightarrow \mathbb{R}^n$ such that the Hamiltonian f\/low of $a_i \circ H_i$ is periodic for each $1\leq i\leq n$.  The map $a$ is often referred to as \emph{global action coordinates} for $\mathcal{H}$.
\end{Definition}

If an integrable system with compact f\/ibers admits global action coordinates, then  $(M,\omega$, $a\circ \mathcal{H})$ is a symplectic toric manifold. Thus, an important question in the study of integrable systems~is

\begin{Question}\label{question}
  Does a given integrable system $(M,\omega,\mathcal{H})$ admit global action coordinates?
\end{Question}

Suppose the f\/ibers of $\mathcal{H}$ are compact and connected; then the regular f\/ibers are Lagrangian $n$-tori and the Arnold--Liouville theorem implies that the restricted map $\mathcal{H}\colon M_{\rm reg} \rightarrow B_{\rm reg}$ is a~f\/iber bundle of compact, connected Lagrangian tori $T^n$ that is locally symplectically equivalent to the trivial system
\begin{gather*}
\mathcal{H} = (x_1, \ldots ,x_n)\colon \ \big(D^n\times T^n, dx\wedge d\theta\big)\rightarrow D^n.
\end{gather*}
Regular integrable systems such as this are called \emph{Lagrangian torus fibrations}. Given a Lag\-ran\-gian torus f\/ibration, there is a natural f\/iber-wise action
\begin{gather*}
T^*B \times_{B} M \rightarrow M
\end{gather*}
def\/ined by $(\alpha_b,m) \mapsto \phi^1_{\alpha_b}(m)$ where $\phi^1_{\alpha_b}$ is the time-1 f\/low of the vector-f\/ield $\omega^{-1}\mathcal{H}^*\alpha$ for some 1-form $\alpha \in \Omega^1(B)$ such that $\alpha(b) = \alpha_b$. The stabilizer of this action is a smooth submanifold $\Lambda\subset T^*B$ that is a full rank lattice in each f\/iber, called the \emph{period lattice} of the f\/ibration. Restricting the projection map $\pi\colon T^*B \rightarrow B$ to this lattice, we obtain a covering space $\pi\colon \Lambda\rightarrow B$ whose topological monodromy is a group homomorphism $m_b\colon  \pi_1(B,b) \rightarrow GL(n,\mathbb{Z})$ after identifying the f\/iber $\Lambda_b$ with $\mathbb{Z}^n$ by a choice of basis (see, e.g., \cite{mein} or \cite{duistermaat} for full details of this construction). The \emph{topological monodromy} of a Lagrangian torus f\/ibration is the topological monodromy of its period lattice, and its non-triviality is the obstruction to the torus f\/iber bundles
\begin{gather*}
  T^*B/\Lambda \rightarrow B\qquad \mbox{and} \qquad M \rightarrow B
\end{gather*}
being principal, and thus to $\mathcal{H}$ admitting global action coordinates. In his seminal 1980 paper \emph{On global action-angle coordinates}, Duistermaat proved

\begin{Theorem}[\cite{duistermaat}]
  A Lagrangian torus fibration $\mathcal{H}\colon M \rightarrow B$ admits global action coordinates if and only if the topological monodromy of the period lattice is trivial.
\end{Theorem}

An integrable system with compact, connected f\/ibers admits global action coordinates only if the corresponding Lagrangian torus f\/ibration $M_{\rm reg} \rightarrow B_{\rm reg}$ admits global action coordinates, so to answer Question~\ref{question} one may begin by studying the topological monodromy of this f\/ibration. On the other hand, there is a subtle interplay between the topological mono\-dromy of \mbox{$M_{\rm reg} \rightarrow B_{\rm reg}$} and the critical f\/ibers of~$\mathcal{H}$. To describe this interplay we need to recall several def\/initions and theorems from~\cite{izosimov,zung1}.

Let $p$ be a critical point of rank $n-k$ of an integrable system $\mathcal{H} = (H_1, \ldots ,H_n)$ on~$(M,\omega)$. Without loss of generality, we may assume that $dH_1,\ldots ,dH_{k} = 0$. The operators $\omega^{-1}d^2H_1$, $\ldots,\omega^{-1}d^2H_{k}$ form a commutative subalgebra~$\mathfrak{h}$ of $\mathfrak{sp}(L^{\perp}/L) \cong \mathfrak{sp}(\mathbb{R},2k)$ where the sub\-space \mbox{$L\subset T_pM$} is the span of the vector f\/ields~$X_{H_{k+1}} ,\ldots , X_{H_n}$ and~$L^{\perp}$ is its symplectic ortho-complement. A~generic critical point of~$\mathcal{H}$ will satisfy the following Morse--Bott type condition for integrable systems:

\begin{Definition}\label{non-degenerate}
  A rank $n-k$ critical point $p$ of an integrable system $\mathcal{H}$ is \emph{non-degenerate} if the subalgebra $\mathfrak{h}$ def\/ined above is a Cartan subalgebra. Similarly, a critical f\/iber of an integrable system is \emph{non-degenerate} if all of its critical points are non-degenerate.
\end{Definition}

   Equivalently, $p$ is non-degenerate if some linear combination of the operators $\omega^{-1}d^2H_1, \dots$, $\omega^{-1}d^2H_{k}$ has $2k$ distinct eigenvalues (and the operators are independent). It was shown in~\cite{williamson}  that a Cartan subalgebra $\mathfrak{h}\subset \mathfrak{sp}(\mathbb{R},2k)$ decomposes as a direct sum of three types of Lie subalgebra, called \emph{elliptic}, \emph{hyperbolic}, and \emph{focus-focus}, and the conjugacy class of a~given Cartan subalgebra~$\mathfrak{h}$ can be described by the three positive integers~$h_e$, $h_h$, and~$h_f$, which denote the number of elliptic, hyperbolic, and focus-focus blocks respectively in this decomposition\footnote{See   Appendix~\ref{appendixA} for a more explicit description of these subalgebras.}.

\begin{Definition}\label{type point}
  The \emph{Williamson type} of a non-degenerate critical point~$p$ for an integrable system $\mathcal{H}$ is the 3-tuple $(h_e,h_h,h_f)$ corresponding to the Cartan subalgebra $\mathfrak{h} \subset \mathfrak{sp}(\mathbb{R},2k)$.
\end{Definition}

The remainder of this paper is primarily concerned with critical points for which the Wil\-liam\-son type is $(h_e,0,h_f)$; such critical points are sometimes called \emph{almost toric}. A~f\/iber is called \emph{almost toric} if all the critical points it contains are almost toric.

It was shown in \cite{zung1} that the Williamson type of a critical point of lowest rank in a non-degenerate critical f\/iber is an invariant of the f\/iber, which allows the def\/inition,

\begin{Definition}\label{type fiber}
  The \emph{Williamson type} of a non-degenerate critical f\/iber $N$ in an integrable system $(M,\omega,\mathcal{H})$ with compact and connected f\/ibers is the Williamson type of any critical point of lowest rank in the f\/iber. The \emph{rank} of $N$ is the rank of these critical points.
\end{Definition}

Assuming once again that all the f\/ibers of an integrable system are compact and connected, we are interested in the foliation of a neighbourhood $\mathcal{U}(N)$ of a critical f\/iber $N$ by the f\/ibers of $\mathcal{H}$ (this is the \emph{Liouville foliation}, although it is def\/ined dif\/ferently in more general settings, see~\cite{zung1}). We consider neighbourhoods of the form $\mathcal{U}(N) = \mathcal{H}^{-1}(D)$ where $D$ is a small disc centred at $\mathcal{H}(N)$. Two foliations are said to be \emph{topologically equivalent} if they are related by a foliation preserving homeomorphism. A~\emph{singularity} of an integrable system is an equivalence class of a Liouville foliation in a neighbourhood of~$N$.

  Let $(H_1,H_2)$ be an integrable system with compact, connected f\/ibers on a symplectic 4-manifold $M$. The (topological equivalence class of~a) Liouville foliation in a neighbourhood of critical f\/iber~$N$ whose only critical points are non-degenerate, rank 0, focus-focus is called a~\emph{stable focus-focus singularity}. The critical f\/iber~$N$ of such a singularity is a disjoint union of $c\geq 1$ critical points $x_1,\ldots,x_c$ and $c$ open annuli~$A_1,\ldots A_c$, $A_i \cong \mathbb{R}\times S^1$, such that each annulus~$A_i$ has $\{x_i,x_{i+1}\}$ as its boundary.

\begin{Theorem}[\cite{zung2}]\label{zung focus focus}
  If $N$ is a stable focus-focus singularity of an integrable system which contains $c\geq 1$ critical points then there is a neighbourhood $\mathcal{U}(N) = \mathcal{H}^{-1}(D^2)$ such that $\mathcal{H}(N) = 0$ is the only critical value in~$D^2$, and the topological monodromy of the torus fibration $\mathcal{U}(N)\setminus N \rightarrow D^2\setminus\{0\}$ is generated by
 \begin{gather*} A = \left(\begin{matrix}
1 & c \\
0 & 1
\end{matrix}\right)
\end{gather*}
$($i.e., up to a choice of basis for $\mathbb{Z}^n$, $m_b(\gamma) = A$ for a representative~$\gamma$ of a generator of $\pi_b(D^2\setminus\{0\}))$.  Furthermore, every two such Liouville foliations are topologically equivalent.
\end{Theorem}

Note that since $m_b$ is a group homomorphism,
\begin{gather*}
m_b\big(\gamma^{-1}\big) = A^{-1} =  \left(\begin{matrix}
1 & -c \\
0 & 1
\end{matrix}\right),
\end{gather*}
but $A^{-1}$ and $A$ are in the same $GL(n,\mathbb{Z})$ conjugacy class.

In particular, a model Liouville foliation $\mathcal{U}(N_c^f) \rightarrow D^2$ can be constructed for any \mbox{$c \geq 1$}~\cite{zung2}. It has been shown that this prototypical model can be used to understand almost toric Liouville foliations in higher degrees of freedom up to topological equivalence, under an additional assumption:

\begin{Definition}[\protect{\cite[Def\/inition~6.3]{zung1}}]\label{topological stability}
A non-degenerate singularity $\mathcal{H}\colon \mathcal{U}(N)\rightarrow D^n$ of an integrable system is called \emph{topologically stable} if the local critical value set of the moment map restricted to $\mathcal{U}(N)$ coincides with the critical value set of the moment map restricted to a small neighbourhood of a singular point of minimal rank in $N$.
\end{Definition}

The following theorem was proven for more general singularities, but we only state it for almost toric f\/ibers.

\begin{Theorem}[\cite{zung1}]\label{zung classification}
  Suppose $N$ is a rank $n-k$, topologically stable, almost toric critical fiber of an integrable system~$\mathcal{H}$ with compact and connected fibers. Then the Liouville foliation of~$\mathcal{U}(N)$ is topologically equivalent to a finite quotient of the product foliation:
  \begin{gather*}
    \big(\big(D^{n-k}\times T^{n-k}\big) \times \mathcal{U}\big(N_{h_e}^e\big) \times  \mathcal{U}\big(N_{c_1}^f\big) \times \cdots \times \mathcal{U}\big(N_{c_{h_f}}^f\big)\big)/G,
  \end{gather*}
  where $G$ is a finite group which acts freely and component-wise on each factor in a foliation preserving manner\footnote{Here $\mathcal{U}(N_{h_e}^e)$ is the model Liouville foliation $D^{2h_e} \rightarrow \mathbb{R}^{h_e}$ given by $(x_1,y_1, \ldots, x_{h_e},y_{h_e}) \mapsto (x_1^2 + y_1^2, \ldots$, $x_{h_e}^2 + y_{h_e}^2)$ as proven by~\cite{eliasson}.}. Moreover, $G$ acts trivially on the elliptic component.
\end{Theorem}

This result was later strengthened by Izosimov,

\begin{Theorem}[\cite{izosimov}]\label{izosimov}
  The Liouville foliation of $\mathcal{U}(N)$ in the preceding theorem is topologically equivalent to the product foliation of
  \begin{gather*}
    \big(D^{n-k}\times T^{n-k}\big) \times \mathcal{U}\big(N_{h_e}^e\big)\times \big(\mathcal{U}\big(N_{c_1}^f\big) \times \cdots \times \mathcal{U}\big(N_{c_{h_f}}^f\big)\big)/G.
  \end{gather*}
\end{Theorem}

\looseness=-1
Given an almost toric critical f\/iber, this theorem implies that the associated Lagrangian torus f\/ibration is topologically equivalent to the associated Lagrangian torus f\/ibration of a~model pro\-duct system. In particular, both f\/ibrations will have the same topological monodromy, since it is a topological invariant of torus f\/ibrations.  The topological mo\-no\-dromy of a~product of Lag\-rangian torus f\/ibrations $\mathcal{H}\times\mathcal{H}'\colon M\times M' \rightarrow B \times B'$ is a the topological mo\-no\-dromy of the lattice $\Lambda_{\mathcal{H}\times \mathcal{H}'} \subset T^*B \times T^*B'$ which decomposes as the product covering space $\Lambda_{\mathcal{H}}\times \Lambda_{\mathcal{H}'} \rightarrow B\times B'$, and the topological mono\-dromy of this covering space decomposes~as
\begin{gather*}
m_b(\gamma) =
\left(\begin{matrix}
m_{b,\mathcal{H}}(\gamma)&0\\
0&m_{b,\mathcal{H}'}(\gamma)
\end{matrix}\right).
\end{gather*}

\begin{Remark}
  In Section~\ref{section 7} we will use Theorem~\ref{izosimov} to conclude that (for the Heisenberg spin system) the Liouville foliation above the critical line is homeomorphic to a~product, and thus the topological monodromy around the critical line decomposes as a~product, as in the preceding discussion.  In fact, one can prove this more directly by using the description in~\cite{zung1} of the group $G$ appearing in Theorem~\ref{zung classification}, and our explicit identif\/ication of the singular f\/ibers with the product~$S^1 \times N^f_3$ (Remark~\ref{identification}) to show that the group~$G$ will vanish, and hence the foliation is homeomorphic to a~product.
\end{Remark}

\section{Symplectic geometry of coadjoint orbits}\label{section 3}

Let $G$ be a compact, connected Lie group with Lie algebra $\mathfrak{g}$ endowed with an Ad-invariant inner product~$\langle\,,\,\rangle$.  After $G$-equivariant identif\/ication of $\mathfrak{g}$ with $\mathfrak{g}^*$ via the inner product, the Kostant--Kirillov--Souriau symplectic
structure on an adjoint orbit $\mathcal{O}_Z$ of $Z\in \mathfrak{g}$ is given by
\begin{gather*}
  \omega_Z([X,Z],[Y,Z]) = \langle Z,[X,Y]\rangle,
\end{gather*}
where $X,Y \in \mathfrak{g}$.
Hamilton's equation for a function $H\colon \mathfrak{g} \rightarrow \mathbb{R}$ can be written
as
\begin{gather*}
  \frac{dZ}{dt} = X_H(Z) =  [\nabla H(Z),Z],
\end{gather*}
where $\nabla H$ is the gradient vector f\/ield def\/ined by the equation $\langle \nabla H(Z), Y\rangle = dH_Z(Y)$ for all $Y\in T_Z\mathcal{O}$. Hence, the Poisson bracket of two functions $H$, $F$ can be conveniently written as
\begin{gather*}
  \{H,F\}_Z = \omega_Z(X_H,X_F) =  \omega_Z([\nabla H, Z],[\nabla F, Z]) = \langle Z,[\nabla H,\nabla F]\rangle.
\end{gather*}
 A direct sum of semisimple Lie algebras $\mathfrak{g}_1\oplus \cdots \oplus \mathfrak{g}_N$ is endowed with direct sum Lie brackets and Killing forms.
 An adjoint orbit in $\mathfrak{g}_1\oplus \cdots \oplus \mathfrak{g}_N$ is a~product of adjoint orbits  $\mathcal{O}_{Z_1}\times \cdots \times \mathcal{O}_{Z_N}$ and the symplectic structure coincides with the direct sum of their respective symplectic structures, $\omega = \omega_1\oplus \cdots \oplus \omega_N$.

The moment map for the adjoint action of $G$ on an orbit in $\mathfrak{g}$ is inclusion of the orbit into $\mathfrak{g}^*$ via the ${\rm Ad}$-equivariant identif\/ication. The moment map for the diagonal adjoint action of $G$ on $\mathcal{O}_{Z_1}\times \cdots \times \mathcal{O}_{Z_N}$ is hence the map $(X_1, \ldots ,X_N) \mapsto \sum X_i$.

\begin{Example}\label{so3}
  Let $G = {\rm SO}(3)$ be the group of rotations of $\mathbb{R}^3$ equipped with the standard basis and inner product. Its Lie algebra is
  \begin{gather*} \mathfrak{so}(3) = \left\{ \left(\begin{matrix}
      0 & x_1 & x_2 \\
      -x_1 & 0 & x_3\\
      -x_2 & -x_3 & 0
    \end{matrix}\right)\colon \, x_1,x_2,x_3\in \mathbb{R}\right\},
  \end{gather*}
  which has the ${\rm Ad}$-invariant inner product $\langle X,Y\rangle = -\frac{1}{2}\operatorname{tr}(XY)$ and Lie bracket $[X,Y] = XY-YX$. The map $\Psi\colon \mathfrak{so}(3) \rightarrow (\mathbb{R}^3,\times)$ given by
  \begin{gather*}
    \left(\begin{matrix}
      0 & x_1 & x_2 \\
      -x_1 & 0 & x_3\\
      -x_2 & -x_3 & 0
    \end{matrix}\right) \mapsto (x_1,x_2,x_3) \in \mathbb{R}^3
  \end{gather*}
  is an isomorphism of Lie algebras, where $(\mathbb{R}^3,\times)$ is the cross-product Lie algebra.  Under this identif\/ication, the adjoint action of ${\rm SO}(3)$ on $\mathfrak{so}(3)$ is identif\/ied with the standard action on~$\mathbb{R}^3$, and the ${\rm Ad}$-invariant inner product is simply the standard Euclidean inner product:
  \begin{gather*}
    \langle X,Y\rangle = x_1y_1 + x_2y_2 + x_3y_3 = -\frac{1}{2}\operatorname{tr}(XY).
  \end{gather*}
  Thus, the adjoint orbits are identif\/ied with concentric spheres in $\mathbb{R}^3$ and the symplectic structure on an adjoint orbit through the point $Z = (z_1,z_2,z_3)$ is
  \begin{gather*}
    \omega_Z([X, Z],[Y, Z]) = \langle Z,[X,Y]\rangle = \det\left(\begin{matrix}
      z_1 & z_2 & z_3 \\
      x_1 & x_2 & x_3\\
      y_1 & y_2 & y_3
    \end{matrix}\right),
  \end{gather*}
  which is precisely the standard symplectic structure on a sphere with radius~$|Z|$.
In what follows, we will consider the adjoint orbit in $\mathfrak{so}(3)\oplus \mathfrak{so}(3) \oplus \mathfrak{so}(3)$ that is a product of three spheres with radius~1, $(S^2\times S^2\times S^2, \omega_{\rm STD}\oplus \omega_{\rm STD}\oplus \omega_{\rm STD})$.
\end{Example}

\section{An integrable Heisenberg spin chain}\label{section 4}

As in Example~\ref{so3}, consider a product of three spheres of radius 1 whose elements are triples $(X,Y,Z)\in M = S^2\times S^2\times S^2$. Def\/ine the Hamiltonians
\begin{gather*}
\begin{split}
&   H(X,Y,Z) = |X+Y+Z| = \sqrt{3+ \langle X,Y\rangle + \langle Y,Z\rangle + \langle Z,X\rangle},
\\
&  I(X,Y,Z)  =  \langle X+Y+Z, e_3\rangle, \qquad \mbox{and} \qquad J(X,Y,Z) = \langle X,[Y,Z]\rangle= \det (X,Y,Z).
\end{split}
\end{gather*}
By ad-invariance of the inner product, the Hamiltonian vector f\/ields of these functions are
\begin{gather*}
  X_H = [ \nabla H,(X,Y,Z)] = \frac{1}{H(X,Y,Z)}([Y+Z,X],[X+Z,Y],[X+Y,Z]),\\
     X_{I} = [\nabla I,(X,Y,Z)] = ([e_3,X],[e_3,Y],[e_3,Z]),\qquad \mbox{and}\\
 X_{J} = [\nabla J,(X,Y,Z)] = ([[Y,Z],X],[[Z,X],Y],[[X,Y],Z]).
\end{gather*}
The f\/low $\varphi_{[v,X]}^t$ of a vector f\/ield $[v,X]$ acts by rotation of the vector $X$ around the axis $v$ with period~$2\pi/|v|$, which we denote as~$R_t^v$,
\begin{gather*}
  \varphi_{[v,X]}^t X = R_t^vX.
\end{gather*}
Thus, the Hamiltonian f\/low of $I$ acts by rotating each sphere around the $e_3$-axis with period $2\pi$,
\begin{gather*}
  \varphi_{X_{I}}^t (X,Y,Z) = \big(R_t^{e_3}X,R_t^{e_3}Y,R_t^{e_3}Z\big).
\end{gather*}
Where def\/ined, the Hamiltonian f\/low of $H$ rotates each sphere around the axis $X+Y+Z$ with period $2\pi$,
\begin{gather*}
  \varphi_{X_{H}}^t(X,Y,Z) = \big(R_t^{v}X,R_t^{v}Y,R_t^{v}Z\big),
\end{gather*}
where $v = (X+Y+Z)/|X+Y+Z|$.
This is perhaps best visualized as rotating the polygon with edges $X$, $Y$, $Z$, $-X-Y-Z$ around the edge $-X-Y-Z$.

\begin{Proposition}\label{commute} $\{H,J\} = \{H,I\} = \{J,I\} = 0$.
\end{Proposition}

\begin{proof} It is a nice exercise to see that this is true based on the geometric description of the Hamiltonians and their f\/lows given above. More algebraically, one can see this using the Lie algebra structure that is present. For example,
\begin{gather*}
\{J,I\}_{(X,Y,Z)}  =  \langle (X,Y,Z), [([Y,Z],[Z,X],[X,Y]), (e_3,e_3,e_3)]\rangle \\
\hphantom{\{J,I\}_{(X,Y,Z)}}{} = \langle [X,[Y,Z]]+[Y,[Z,X]] + [Z,[X,Y]], e_3\rangle = 0
\end{gather*}
by the Jacobi identity and ad-invariance of our inner product. A similar calculation shows that $\{H,I\} = \{H,J\} =0$.
\end{proof}

\begin{Remark}\label{lagrangian rp3}
The Hamiltonian f\/low of $J$ is less straightforward to describe, but we can say something about how it acts on the submanifold
\begin{gather*}
L = \{(X,Y,Z) \in M \colon  X + Y + Z = 0 \} = H^{-1}(0).
\end{gather*}
There the vectors $X$, $Y$, and $Z$ are coplanar and the vector f\/ield
\begin{gather*}
X_{J}(X,Y,Z)= ([[X,Y],X],[[Y,Z],Y],[[Z,X],Z]).
\end{gather*}
The f\/low of this vector f\/ield acts on $L$ by rotation of each vector $X$, $Y$ and $Z$ around the axis $\overline{n} = [X,Y] = [Y,Z] = [Z,X]$ with constant period: since each component of~$X_{J}$ is tangent to the plane spanned by~$X$,~$Y$, and~$Z$, and the f\/low of $X_{J}$ preserves~$L$, the vector $[X,Y]$ is preserved by~$X_{J}$, so
\begin{gather*}
X_J = ([\overline{n},X],[\overline{n},Y],[\overline{n},Z]).
\end{gather*}

In \cite{oakley} it was shown that the f\/iber $L$ is a non-displaceable Lagrangian submanifold of~$(S^2\times S^2\times S^2, \omega_{\rm STD}\oplus \omega_{\rm STD} \oplus \omega_{\rm STD})$. To see that it is an embedded Lagrangian~$\mathbb{R} P^3$, observe that it is the zero level set for the moment map of the diagonal ${\rm SO}(3)$-action, $(X,Y,Z) \mapsto X+Y+Z$, and the diagonal action of~${\rm SO}(3)$ is free and transitive.
\end{Remark}

\begin{Remark}\label{symmetry}
The Hamiltonian $H$ has several important symmetries. First, note that $H$ is invariant under the diagonal Hamiltonian ${\rm SO}(3)$-action, whose moment map is
\begin{gather*}
  S^2\times S^2 \times S^2 \xrightarrow{(X,Y,Z)\mapsto X+Y+Z} \mathfrak{so}(3)\cong \mathfrak{so}(3)^*.
\end{gather*}
Thus $H$ is non-commutatively integrable, or super-integrable. Indeed, there are three independent choices for our second integral
\begin{gather*}
  I_v(X,Y,Z) = \langle X+Y+Z,v\rangle,
\end{gather*}
which do not pairwise Poisson commute.

In addition, all three Hamiltonians, $H$, $I$, $J$ have a symplectic $\mathbb{Z}_3$-symmetry by cyclic permutations $(X,Y,Z) \mapsto (Z,X,Y)$. As we will see, this symmetry of the system by a~group of order~3 is also a symmetry of the rank~1 focus-focus critical f\/ibers, and thus the system's topological monodromy contains the number~3 (see Theorem~\ref{zung focus focus}).
\end{Remark}


\section{Image of the moment map}\label{section 5}

In this section we give a complete description the image, critical set, and  critical values for the moment map $\mathcal{H}= (H,I,J)$.

The image of the moment map has two ref\/lective symmetries coming from the fact that $J(X,Y,Z) = -J(X,Z,Y)$ and  $I(-X,-Y,-Z) = - I(X,Y,Z)$. It is obvious that $|I| \leq H$, with equality when $X+ Y + Z\in \operatorname{span}(e_3)$, and that $0\leq H \leq 3$ with equality when $X = Y = Z$.

Observe that
\begin{gather*}
  H = \sqrt{3 +2(a+b+c)} \qquad \mbox{and} \qquad |J| = \sqrt{1 + 2abc -\big(a^2+b^2+c^2\big)},
\end{gather*}
where $a = \langle X,Y\rangle$, $b = \langle Y,Z\rangle$ and $c = \langle Z,X\rangle$ (the second formula is the volume of a parallelpiped). If we maximize $|J|$ with the constraint $H = \operatorname{const}$, then we must have $a = b =c$ (the interior angles between the three vectors are the same). Using this we can deduce that
\begin{gather*}
|J| \leq \sqrt{1 + 2\left(\frac{H^2-3}{6}\right)^3 - 3\left(\frac{H^2-3}{6}\right)^2},
\end{gather*}
 and equality is achieved when $a=b=c$.

\begin{Proposition}\label{moment map image}
  The image of the moment map $\mathcal{H}$ is the set of points $(r,s,t)\in \mathbb{R}^3$ that satisfy the equations
  \begin{gather*}
    |t| \leq \sqrt{1 + 2\left(\frac{r^2-3}{6}\right)^3 - 3\left(\frac{r^2-3}{6}\right)^2}, \qquad |s| \leq r, \qquad \mbox{and} \qquad 0 \leq  r \leq 3.
  \end{gather*}
\end{Proposition}

\begin{proof} Observe that the level sets of $H$ are all connected (see the proof of Theorem~\ref{level sets} in Section~\ref{section 6}). For a given tuple $(r,s,t)\in \mathbb{R}^3$ that satisf\/ies the inequalities of~\eqref{region}, consider the restriction of the map $J$ to the connected set $H^{-1}(r)$. Since $0 \leq r \leq 3$, we can construct a tuple $(X,Y,Z)\in H^{­1}(r)$ such that
\begin{gather*}
  \langle X,Y\rangle = \langle Y,Z\rangle = \langle Z,X\rangle = \frac{r^2 -­3}{6}.
\end{gather*}
$J(X,Y,Z)$ is maximal or minimal depending on the orientation of the basis $\{X,Y,Z\}$ (if $r = 0$ or $3$ then $\{X,Y,Z\}$ does not form a~basis and $J(X,Y,Z)=0$ is both maximal and minimal). It follows by the Intermediate Value Theorem that there is a tuple $(X,Y,Z) \in S^2 \times S^2 \times S^2$ such that $H(X,Y,Z) = r$ and $J(X,Y,Z) = t$. Using the Intermediate Value Theorem again, it is easy to see that there exists a $\theta\in [0,2\pi]$ such that  $I(R^{e_3}_{\theta}X,R^{e_3}_{\theta}Y,R^{e_3}_{\theta}Z) = s$. Since $H$ and $J$ are invariant under diagonal rotations, we are done.
\end{proof}

Next, we turn our attention to the critical set for the system. The critical set consists of several subsets:
\begin{enumerate}\itemsep=0pt
\item The sets where $H$ is critical:
\begin{enumerate}\itemsep=0pt
\item the embedded ${\rm SO}(3) \cong H^{-1}(0)$, which lies over the vertex $(0,0,0)$,
\item the three embedded spheres $S_1 = \{(-X,X,X)\}$, $S_2 = \{(X,-X,X)\}$, and $S_3 = \{(X,X,-X)\},$ which lie over the critical line ($H=1$ and $J = 0$), and
\item the diagonally embedded sphere $S_4 = \{(X,X,X)\}$, which lies over the edge $H=3$.
\end{enumerate}
\item The set $C_1 = \{(X,Y,Z)\colon  X + Y + Z \in \operatorname{Span}(e_3)\}$ where $dI$ is proportional to~$dH$. This contains the critical set of $I$ and is mapped to two opposite faces of the moment map image.
\item The set $C_ 2 = \{(X,Y,Z)\colon \langle X,Y\rangle =  \langle Y,Z\rangle = \langle Z,X\rangle \}$ where $dJ$ is proportional to  $dH$. Note that $S_4, H^{-1}(0) \subset C_2$. This set maps to the other two opposite faces of the moment map image. Points in $C_1 \cap C_2$ map to edges of the moment map image.
\end{enumerate}

\begin{Proposition}\label{critical} The critical set for $\mathcal{H}$ is
\begin{gather*}C = C_1 \cup C_2 \cup S_1 \cup S_2 \cup S_3.\end{gather*}
The set $\mathcal{H}(C_1 \cup C_2)$ is the boundary of the image $\mathcal{H}(M)$ and the set $ \mathcal{H}(S_1) = \mathcal{H}(S_2)=\mathcal{H}(S_3)$ is the line segment $\mathcal{H}(M) \cap \{ H = 1,J = 0\}$ $($see Fig.~{\rm \ref{figure 1})}.
\end{Proposition}

In the terminology of integrable systems, the set of critical values is the system's `bifurcation diagram'. This description of the bifurcation diagram will be of use to us in Section~\ref{section 7}. Note that the three critical spheres $S_1$, $S_2$, $S_3$ are permuted by the system's $\mathbb{Z}_3$-symmetry (cf.\ Remark~\ref{symmetry}).

\begin{proof} Throughout we use the fact that $df = 0$ if and only if $X_f= 0$ (when $f$ is smooth).
\begin{enumerate}\itemsep=0pt
\item Since $0$ is the global minimum for $H$, the set $H^{-1}(0)$ is critical.  To f\/ind the other critical sets of $H$, observe that $X_{H} = 0$ if and only if $[X+Y+Z,X] = 0$,  $[X+Y+Z,Y] = 0$, and $[X+Y+Z,Z] = 0$. If $X + Y + Z \neq 0$, this occurs if and only if $X+Y+Z$ is is contained in the lines spanned by~$X$, $Y$, and~$Z$, so this happens if and only if~$X$, $Y$, and~$Z$ are collinear. This entails cases~(1b) and~(1c).

\item Observe that $X_{I} = \alpha X_{H}$ for some $\alpha\neq 0$ if and only if
\begin{gather*}
    [e_3,X]  = \alpha[X+Y+Z,X],\qquad     [e_3,Y]  = \alpha[X+Y+Z,Y] , \qquad \mbox{and}\\
    [e_3,Z]  = \alpha[X+Y+Z,Z],
\end{gather*}
which is true if and only if $X+Y+Z \in \operatorname{span}(e_3)$.
\item Observe that $X_{J} = \alpha X_{H}$ if and only if
\begin{gather*}
    [[Y,Z],X] = \alpha[Y+Z,X],\qquad
    [[Z,X],Y]  = \alpha [X +Z,Y], \qquad \mbox{and}\\
    [[X,Y],Z]  = \alpha[X+Y,Z]
\end{gather*}
for some $\alpha$. If $\alpha = 0$ then the 3-tuple~$X$,~$Y$,~$Z$ forms an oriented or anti-oriented orthonormal frame, or~$X$,~$Y$ and~$Z$ are collinear.  If $\alpha \neq 0$ then
\begin{enumerate}\itemsep=0pt
  \item $[Y,Z]$, $Y+Z$, and $X$ are coplanar,
  \item $[Z,X]$, $X+Z$, and $Y$ are coplanar, and
  \item $[X,Y]$, $X+Y$, and $Z$ are coplanar.
\end{enumerate}
Since $X$, $Y$ and $Z$ all have the same length, the three items listed above are true if and only if $X$, $Y$, $Z$ are collinear, or  $\langle X,Y\rangle = \langle Y,Z \rangle = \langle Z,X\rangle$ (this fact is a straightforward exercise in Euclidean geometry).
\end{enumerate}

Finally, suppose that  $\alpha X_{I}+\beta X_{H} + \gamma X_{J} = 0$ for $\alpha,\beta, \gamma \in \mathbb{R}$ not all zero. Then
\begin{gather*}
    \alpha [e_3,X] + \beta [Y + Z,X] + \gamma[[Y,Z],X] = 0,\\
    \alpha [e_3,Y] + \beta [X + Z,Y] + \gamma[[Z,X],Y] = 0, \qquad \mbox{and}\\
    \alpha [e_3,Z] + \beta [X + Y,Z] + \gamma[[X,Y],Z] = 0.
\end{gather*}
Adding these equations together we obtain
\begin{gather*}
  \alpha [e_3,X+Y+Z] +  \gamma\left([[X,Y],Z] + [[Y,Z],X] + [[Z,X],Y]\right) = 0,
\end{gather*}
which by the Jacobi identity reduces to $\alpha [e_3,X+Y+Z] = 0$. If $[e_3,X+Y+Z]=0$, then $(X,Y,Z) \in C_1$. If  $\alpha = 0$, then $(X,Y,Z)$ is in the set  $S_1\cup S_2\cup S_3 \cup C_2$.
\end{proof}

Combining Propositions~\ref{moment map image} and \ref{critical}, we have Theorem~\ref{image of moment map}. Since we now know that our Hamiltonians are independent on an open dense subset of $M$, we can conclude:

\begin{Corollary}
$\mathcal{H}$ is a completely integrable system. In particular, the regular level sets of $\mathcal{H}$ are homeomorphic to a disjoint union of finitely many $3$-tori.
\end{Corollary}

\begin{Corollary} The set of regular values is homotopy-equivalent to~$S^1$.
\end{Corollary}

In the next section, we will see that the regular level sets are connected and in Section~\ref{section 7} we will describe the structure of the associated Lagrangian torus f\/ibration.

\section{Connectedness of regular level sets}\label{section 6}

In this section we prove

\begin{Theorem}\label{level sets}
  The regular fibers of the map $\mathcal{H} = (H,I,J)$ are all connected.
\end{Theorem}

\begin{proof}
The proof will have two parts. For $H \neq 1$ we can make a general argument and for $H = 1$ we will apply Ehresmann's theorem.

Pick a regular value $(r,s,t)$  with $r\neq 1$. Since $H$ generates a Hamiltonian $S^1$-action on $M\setminus H^{-1}(0)$, it is a Morse--Bott function such that all critical sets have even index \cite{atiyah}. Thus, the level sets of $H$ are connected. Since the regular level sets of $H$ are compact and connected, the symplectic reductions $M_r \equiv H^{-1}(r)/S^1$ are all compact, connected symplectic manifolds. Since $s$ is a regular value of the reduced Hamiltonian~$\tilde I$, and $\tilde I$ generates a free $S^1$-action, we can reduce once more to obtain the compact and connected manifold $M_{r,s} = \tilde I^{-1}(s)/S^1$.

The image of the twice-reduced Hamiltonian $\tilde J$ on $M_{r,s}$ is a line segment with the only critical values being the maximum and minimum. Recall from Proposition~\ref{critical} that in the unreduced manifold, $dJ$ is dependent on~$dH$ and~$dI$ at a point $(X,Y,Z)\in H^{-1}(r)\cap I^{-1}(s)$ if and only if~$dJ$ is proportional to~$dH$, and this occurs if and only if~$\langle X,Y \rangle = \langle Y,Z \rangle = \langle Z,X\rangle$. It is easy to see geometrically that the set of all 3-tuples of unit length vectors which satisfy the conditions
\begin{enumerate}\itemsep=0pt
  \item[1)] $  \langle X,Y \rangle  = \langle Y,Z \rangle = \langle Z,X\rangle$,
  \item[2)] $ |X+Y+Z| = r $, and
  \item[3)] $ \langle X+Y+Z,e_3\rangle  = s$.
\end{enumerate}
is two orbits of the Hamiltonian $T^2$-action generated by $(H,I)$, which are distinguished by whether the  basis $\{X,Y,Z\}$ is positively or negatively oriented, corresponding to being the maximum or minimum set for~$J$ on $H^{-1}(r)\cap I^{-1}(s)$. Thus there are two critical points of~$\tilde J$ on~$M_{r,s}$ and the regular f\/ibers $\tilde J^{-1}(s)$ are all connected. This implies that the regular f\/iber $\mathcal{H}^{-1}(r,s,t)\subset M$ is connected.

Now consider a regular value $(1,s,t)$. Since the map $\mathcal{H}$ is proper, Ehresmann's theorem implies that the f\/iber $\mathcal{H}^{-1}(1,s,t)$ is dif\/feomorphic to a nearby f\/iber $\mathcal{H}^{-1}(r,s,t)$ with $r\neq 1$, which we have just shown is connected.
\end{proof}

\begin{Remark} The image of the invariant Lagrangian $L = H^{-1}(0)$ in the reduction at 0 by $I$ is a~Lagrangian~$S^2$.
 \end{Remark}

\begin{Remark} In \cite{prv} it is shown that if $(M^4,\omega,(H_1,H_2))$ is a non-degenerate completely integrable system (cf.\ Def\/inition~\ref{non-degenerate}) with two degrees of freedom, whose set of critical values has no vertical tangent lines, then the system has connected f\/ibers. After rotating the moment map image, and checking that the boundary of~$\mathcal{H}(M)$ consists generically (almost all values of~$H$ or~$I$) of non-degenerate elliptic critical values one could deduce connectedness of almost all the f\/ibers of~$\mathcal{H}$ or~$I$ by applying this theorem to the reduced systems, then applying Ehresmann to the remaining f\/ibers. For example, the systems obtained by reducing at $H=0$ or~$\pm1$ will have degenerate critical points (the Lagrangian $S^2$ of the previous remark and image of intersections $S_i \cap C_1$, respectively), so the theorem of~\cite{prv} cannot be applied directly.
\end{Remark}

\section{Topological monodromy around the critical line}\label{section 7}

As we have seen, the regular f\/ibers of the moment map $\mathcal{H} = (H,I,J)$ for the Heisenberg spin chain are connected and, because of the critical line, the set of regular values $B_{\rm reg}$ is homotopy equivalent to $S^1$. Thus, the associated Lagrangian torus f\/ibration $M_{\rm reg} \rightarrow B_{\rm reg}$ of the Heisenberg spin chain may have non-trivial topological monodromy. This section uses our understanding of the critical set and critical values of $\mathcal{H}$ from Section~\ref{section 5}, together with a non-degeneracy computation that is relegated to Appendix~\ref{appendixA}, to deduce the topological monodromy directly from Theorem~\ref{izosimov}. Recall

\medskip

\noindent\textbf{Theorem 1.1.} \emph{The topological monodromy of the Lagrangian torus fibration associated to the integrable system $(H,I,J)$ is generated by the matrix
\begin{gather*}
A =
\left(\begin{matrix}
1&0&0\\
0&1&3\\
0&0&1
\end{matrix}\right).
\end{gather*}
Thus the system does not admit global action coordinates.}

\begin{proof}
  Consider critical the f\/ibers $N(s) = \mathcal{H}^{-1}(1,s,0)$ for $-1<s<1$.  The condition $J = 0$ implies that the unit-length vectors~$X$,~$Y$, $Z$ are coplanar and the condition $H= 1$ implies that they form three sides of the parallelogram $X$, $Y$, $Z$, $-X-Y-Z$. The set of all such vectors is compact and connected (see the remark following the proof).  By Proposition~\ref{critical}, the set of critical points in the f\/iber~$N(s)$ contains three connected components:
  \begin{gather*}N(s) \cap S_1 = \left\{(-X,X,X)\in S^2 \times S^2 \times S^2\colon \, \langle X, e_3\rangle = s\right\},\\
N(s) \cap S_2 = \left\{(X,-X,X)\in S^2 \times S^2 \times S^2\colon \, \langle X, e_3\rangle = s\right\},\\
N(s) \cap S_3 = \left\{(X,X,-X)\in S^2 \times S^2 \times S^2\colon \, \langle X, e_3\rangle = s\right\},
\end{gather*}
  each homeomorphic to $S^1$, and they are permuted by the system's $\mathbb{Z}_3$ symmetry (cf.\ Remark~\ref{symmetry}). Each critical f\/iber is topologically stable (Def\/inition~\ref{topological stability}) since the system's $\mathbb{Z}_3\times S^1$-symmetry acts transitively on the critical set in each f\/iber~$N(s)$, and preserves the map~$\mathcal{H}$.

By Proposition~\ref{appendix proposition} the critical points in $N(s)$ are rank 1 non-degenerate, with Williamson type~$(0,0,1)$. Thus by Theorem~\ref{izosimov}, the Liouville foliation of $\mathcal{U}(N(s))$ is topologically equivalent to the product foliation of $(D^1 \times S^1) \times \mathcal{U}(N_3^f)$. As noted in Section~\ref{section 2}, this implies that the topological monodromy of our system is the same as the product foliation of $(D^1\times S^1)\times \mathcal{U}(N_3^f)$, which decomposes into blocks for each component of the direct sum,
\begin{gather*}
A =
\left(\begin{matrix}
1&0&0\\
0& a_{22}& a_{23}\\
0& a_{32} & a_{33}\\
\end{matrix}\right).
\end{gather*}
By Proposition~\ref{zung focus focus}, the bottom right minor is
\begin{gather*}
\left(\begin{matrix}
1&3\\
0&1
\end{matrix}\right).\tag*{\qed}
\end{gather*}
\renewcommand{\qed}{}
\end{proof}

Since the Hamiltonian f\/lows of $H$ and $I$ are periodic, the 1-forms $dH$ and $dI$ are global sections of the period lattice $\Lambda \rightarrow B_{\rm reg}$. Thus if $\{dH(b),dI(b),\alpha_b\}$ is a $\mathbb{Z}$-basis for $\Lambda_b$, and $\gamma$ is a choice of generator for $\pi_1(B_{\rm reg},b)$, then we must have
\begin{gather*}
m_b(\gamma) =
\left(\begin{matrix}
1&0&*\\
0&1&*\\
0&0&\pm 1
\end{matrix}\right).
\end{gather*}
Theorem 1.1 tells us that
\begin{gather*}
m_b(\gamma) =
\left(\begin{array}{ccc}
1&0&0\\
0&1&3\\
0&0&1\\
\end{array}\right)
\end{gather*}
up to a choice of basis f\/ixing $dH(b)$ and $dI(b)$.

\begin{Remark}\label{identification}
  It is easy to see how the f\/iber $N(s)$ is homeomorphic to $S^1 \times N_3^f$, where $N_3^f$ is the focus-focus f\/iber with three critical points as described in Section~\ref{section 2}. First note that the set
  \begin{gather*}
    S = \{ (X,Y,Z) \in S^2 \times S^2 \times S^2 \colon \, X + Y + Z = \beta e_1+se_3, \beta \in \mathbb{R}^+\}
  \end{gather*}
  is a slice for the Hamiltonian $S^1$-action generated by $I$ in a neighbourhood of $N(s)$. Since this $S^1$-action is free in a neighbourhood of $N(s)$, the f\/iber $N(s)$ is homeomorphic to $S^1 \times \left(N(s) \cap S\right)$. The set $N(s) \cap S$ consists of: the three annuli,
  \begin{gather*}
  \{(X, \beta e_1+se_3,-X) \in S|\, X \in S^2\setminus\{ \pm( \beta e_1+se_3)\}\},\\
\{( \beta e_1+se_3,X,-X) \in S|\, X \in S^2\setminus\{ \pm( \beta e_1+se_3)\}\},\\
\{(X,-X, \beta e_1+se_3) \in S|\, X \in S^2\setminus\{ \pm( \beta e_1+se_3)\}\}
  \end{gather*}
  together with the three points,
  \begin{gather*}( -\beta e_1-se_3, \beta e_1+se_3,\beta e_1+se_3),\qquad
( \beta e_1+se_3, -\beta e_1-se_3, \beta e_1+se_3),\\
( \beta e_1+se_3, \beta e_1+se_3,- \beta e_1-se_3).
  \end{gather*}
\end{Remark}

\begin{Remark}\label{polygon} The fact that this system has non-trivial monodromy should be unsurprising for the following reason: the topology of the $H$-level sets changes as you pass through the critical value~1. This can be seen directly with Morse theory, but there is also a natural  interpretation in terms of the topology of polygon spaces (as introduced by~\cite{kapmill}). There is a natural dif\/feomorphism of the level set $H^{-1}(r)$ with the manifold $M(1,1,1,r)$ of closed 4-gons in~$\mathbb{R}^3$ with side lengths~$1$,~$1$,~$1,$ and~$r$. When $r \neq 1$, it has been observed by Knutson, Hausmann~\cite{hausknut1}, and Kapovitch and Millson~\cite{kapmill} that the quotient $M(1,1,1,r)/{\rm SO}(3)$ is homeomorphic to~$S^2$, and that the quotient map $\pi\colon M(1,1,1,r) \rightarrow S^2$ is a principal ${\rm SO}(3)$-bundle. Further, the cha\-rac\-teristic classes of this principal ${\rm SO}(3)$-bundle were described by Knutson and Hausmann in their paper~\cite{hausknut2}. Their result says that for $0< r< 1$ the bundle is trivial, whereas for $1< r< 3$ the bundle is non-trivial. It is a fun exercise to check that ${\rm SO}(3)\times S^2$ and the total space of the non-trivial principal ${\rm SO}(3)$-bundle over $S^2$ are not homeomorphic\footnote{Hint: compute their fundamental groups.}.

It was observed in \cite{cushman} that such a change in the topology of the level set $H^{-1}(r)$ as $r$ passes through an interior critical value indicates that there \emph{must} be non-trivial monodromy around the associated critical f\/ibres, since this forces the pullback of the torus bundle to any circle around the critical line to be non-trivial.
\end{Remark}

\begin{Remark}\label{chern} Since $B_{\rm reg}$ is homotopy equivalent to $S^1$, the Chern class of the torus f\/ibration $M_{\rm reg} \rightarrow B_{\rm reg}$ is trivial, so there exists a global section $\sigma\colon B_{\rm reg} \rightarrow M_{\rm reg}$. Since the Lagrangian Chern class \cite{duistermaat} also vanishes, $\sigma$ can be chosen to be Lagrangian, so the map
\begin{gather*}
  \Psi\colon \  T^*(B_{\rm reg})/\Lambda \times_{B_{\rm reg}}M_{\rm reg} \rightarrow M_{\rm reg}
\end{gather*}
gives a symplectomorphism $\alpha\mapsto \Psi (\alpha,\sigma(\pi(\alpha)) )$ which is an isomorphism of the Lagrangian torus f\/ibrations $(T^*(B_{\rm reg})/\Lambda,d\lambda) \rightarrow B_{\rm reg}$ and $(M_{\rm reg},\omega) \rightarrow B_{\rm reg}$ (see~\cite{duistermaat, mein} for details).
\end{Remark}

\begin{figure}[t]
        \centering
\includegraphics{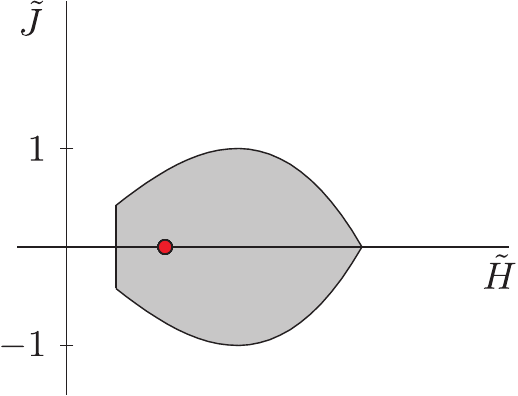}
                \caption{Reduced system on $M_s$ for $-1< s< 1$, $t\neq 0$.}
                \label{figure 3}
\end{figure}

\begin{Remark}
  For $-1< s < 1$, the reduced system on $M_s = I^{-1}(s)/S^1$ has a stable focus-focus critical f\/iber with three critical points. The manifold $M_s$ is the blow-up $Bl_3\mathbb{C} P^2$, and the reduced moment map image has 3 vertices (see Fig.~\ref{figure 3}). A quick comparison with the list of almost toric systems in \cite{leungsym} shows that this checks out.
\end{Remark}

\section{Further directions}

\looseness=-1
In quantum mechanics, the spin-$S$ isotropic Heisenberg spin chain on the lattice $L = \mathbb{Z}$ is an operator $\mathcal{K}$ on the tensor product $V = \bigotimes_{L} \mathbb{C}^{2S+1}$ of irreducible $\mathfrak{su}(2)$ representations, def\/ined as
\begin{gather*}
  \mathcal{K} = \sum_{i\in L} \langle X_i , X_{i+1} \rangle = \sum_{i\in L} \big( X_i^1\circ X_{i+1}^1 + X_i^2\circ X_{i+1}^2+X_i^3\circ X_{i+1}^3\big),
\end{gather*}
where $X_i^j$ are the Pauli matrices
\begin{gather*}
X^1_i  =  \frac{1}{S} \left(\begin{matrix}
0 & -1 \\
1 & 0
\end{matrix}\right),\qquad X^2_i  =  \frac{1}{S} \left(\begin{matrix}
0 & i \\
i & 0
\end{matrix}\right),\qquad X^3_i  =  \frac{1}{S} \left(\begin{matrix}
-i & 0 \\
0 & i \\
\end{matrix}\right)
\end{gather*}
acting on $V$ by the irreducible representation of $\mathfrak{su}(2)$ in the $i$th factor.  This operator models `nearest-neighbour' spin coupling on the lattice. One may consider the inf\/inite chain $L = \mathbb{Z}$ or a chain with boundary conditions, $L = \mathbb{Z}_N$.  For a chain with boundary conditions, the large-$S$ limit of this system is the Hamiltonian
\begin{gather*}
  K = \sum_{i\in L} \langle X_i ,X_{i+1}\rangle
\end{gather*}
on the symplectic manifold  $M_N = (S^2 \times \cdots \times S^2, \omega_{\rm STD} \oplus \cdots \oplus \omega_{\rm STD})$ with elements $(X_1, \ldots ,X_N)$. In the case $N=3$, this is related to the Hamiltonian $H$ studied in this paper by the identity
\begin{gather*}
  3-H^2 = K
\end{gather*}
(of course, the analogous identity fails for $N>3$). In \cite{matthieu}, the algebraic structure of the spin-$\frac{1}{2}$ representation was used to f\/ind a recursive formula for the operators commuting with $\mathcal{K}$ for any lattice $L$ and $S = \frac{1}{2}$. Inspired by this, one might ask if there is a similar combinatorial pattern generating  integrals of the classical Hamiltonians $K$ or $H$ on $M_N$, where we might def\/ine $H_N$ for $N>3$ as the Hamiltonian describing `complete graph' coupling
\begin{gather*}
  H_N = \sum_{i,j\in L,\, i\neq j} \langle X_i ,X_{j}\rangle.
\end{gather*}
\begin{Question}
  How is the structure of the Lie algebra $\mathfrak{so}(3)\cong \mathfrak{su}(2)$ reflected in the conservation laws of the Hamiltonians $H_N$ on $M_N$.
\end{Question}

For example, for arbitrary $N$ one can check using the Jacobi identity and ${\rm ad}$-invariance of~$\langle\,,\,\rangle$, as in the proof of Proposition~\ref{commute} that $H_N$ has integrals
\begin{gather*}
  I_v = \sum_{i\in L}\langle X_i,v\rangle, \quad v\in \mathbb{R}^3,
\qquad
  J =\sum_{i\in L} \langle X_i,[X_{i+1},X_{i+2}]\rangle
\end{gather*}
with $\{H_N,J\} = \{H_N,I_v\} = \{J,I_v\} = 0$, and one could ask if further independent integrals exist. As in Remark~\ref{polygon}, the level sets $H_N^{-1}(h)$ are dif\/feomorphic to the space $M(1,\ldots,1,h)$ of $N$-gons in $\mathbb{R}^3$ with corresponding edge lengths, and come with a~quotient map to the associated polygon space $M(1,\ldots,1,h)/{\rm SO}(3)$.  As~$h$ passes through critical values (which are even or odd integers depending on the parity of~$N$) the topology of this bundle will change, as described in~\cite{hausknut2}, so if~$H_N$ is completely integrable we anticipate non-trivial topological monodromy, provided these critical values lie in the interior of the moment map image.

\appendix

\section{Non-degeneracy computation}\label{appendixA}

In this appendix we prove

\begin{Proposition}\label{appendix proposition}
  The critical points of $\mathcal{H} = (H,I,J)$ which lie above the critical line $\{ (1,s,0)\colon$ ${-}1<s<1\}$ are non-degenerate and have Williamson type $(0,0,1)$ $($see Definitions~{\rm \ref{non-degenerate}} and~{\rm \ref{type fiber})}.
\end{Proposition}

As we have seen, this set of critical points is precisely
\begin{gather*}
\left\{(X,Y,Z)\in S_1\cup S_2\cup S_3 \colon\, I(X,Y,Z) \neq \pm 1\right\}.
\end{gather*}
On this set $dJ = dH = 0$ and $dI \neq 0$, so these points are all rank 1. To show the non-degeneracy of these critical points, we need to show that for each such $p$, the operators $A_J(p) = \omega^{-1}d^2J(p),$ and $A_H(p) = \omega^{-1}d^2H(p)$ span a Cartan subalgebra (see Section~\ref{section 2}). Equivalently, we must f\/ind a linear combination of the operators~$A_J(p)$ and~$A_H(p)$ which has~4 distinct eigenvalues. Note that it is suf\/f\/icient to check non-degeneracy of a single critical point in each f\/iber $\mathcal{H}^{-1}(1,s,0)$ because each connected component of the critical set is an orbit of the Hamiltonian f\/low of~$X_I$, and the symplectic $\mathbb{Z}_3$-action generated by the permutation $\sigma(X,Y,Z) = (Z,X,Y)$ preserves $\mathcal{H}$ and permutes these connected components transitively.

According to the classif\/ication of~\cite{williamson}, the Williamson type of~$p$ will then be determined by the form of these eigenvalues:  in $\mathfrak{sp}(\mathbb{R},4)$ there are four conjugacy classes of Cartan subalgebras corresponding to four possible combinations of eigenvalues for a generic element:
\begin{enumerate}\itemsep=0pt
\item[1)] elliptic-elliptic: $\pm iA$, $\pm iB$,
\item[2)] elliptic-hyperbolic: $\pm A$, $\pm iB$,
\item[3)] hyperbolic-hyperbolic: $\pm A$, $\pm B$, and
\item[4)] focus-focus: $A\pm iB$, $-A\pm iB$.
\end{enumerate}

In Darboux coordinates the operator $A_H(p)$ is equal to the linearization of the Hamiltonian vector f\/ield $X_H$ at~$p$, since
\begin{gather*}
\frac{\partial X_H^i}{\partial x^j }  = \frac{\partial}{\partial x^j} \left( \omega^{ik}\frac{\partial f}{\partial x^k}\right) = \omega^{ik}\frac{\partial^2 f}{\partial x^j\partial x^k} = \big(\omega^{-1}d^2H\big)^i_j.
\end{gather*}
Consider the cylindrical coordinates $(\theta,z) \in (-\pi/2,3\pi/2)\times(-1,1)$ with symplectic form $d\theta_i\wedge dz_i$. The map $\phi\colon (-\pi/2,3\pi/2) \times (-1,1) \rightarrow S^2$ given by
\begin{gather*}
\phi(\theta,z) = \big(\big(1-z^2\big)^{1/2}\cos(\theta),\big(1-z^2\big)^{1/2}\sin(\theta),z\big)
\end{gather*}
is a symplectomorphism. In cylindrical coordinates $(\theta_1,z_1,\theta_2,z_2,\theta_3,z_3)$,  the Hamiltonians are
\begin{gather*}
\hat H   = \left(\sum_j \big(1-z_j^2\big)^{1/2}\cos(\theta_j)\right)^2+\left(\sum_j \big(1-z_j^2\big)^{1/2}\sin(\theta_j)\right)^2 + \left(\sum_j z_j\right)^2,\\
\hat J   = \sum_{j=1,2,3}z_j\big(1-z_{j+1}^2\big)^{1/2}\big(1-z_{j-1}^2\big)^{1/2}\sin(\theta_{j-1}-\theta_{j+1}), \qquad
\hat I  = z_1 + z_2 + z_3
\end{gather*}
(where we have pulled back $(H^2-3)/2$ instead of $H$ for computational convenience). Hamilton's equations tell us that
\begin{align*}
    X_{f} = \frac{\partial f}{\partial z_i}\frac{\partial}{\partial \theta_i}
    - \frac{\partial f}{\partial \theta_i}\frac{\partial}{\partial z_i}.
\end{align*}
The linearization of $X_{f}$ at a f\/ixed point of the f\/low of $f$ is then
\begin{gather*} A_{f} =  \left(\begin{matrix}
\left(\displaystyle\frac{-\partial^2f}{\partial z_k\partial \theta_i}\right)_{ik} & \left(\displaystyle\frac{-\partial^2f}{\partial \theta_k\partial \theta_i}\right)_{ik}\vspace{1mm}\\
\left(\displaystyle\frac{\partial^2f}{\partial z_k\partial z_i}\right)_{ik} & \left(\displaystyle\frac{\partial^2f}{\partial \theta_k\partial z_i}\right)_{ik}
\end{matrix}\right).
\end{gather*}
In order to check non­-degeneracy, one therefore computes the partial derivatives:
\begin{gather*}
\frac{\partial^2 \hat H}{\partial z_k \partial z_i}   =   \begin{cases}
 \displaystyle -2\big(1-z_i^2\big)^{-3/2}\left(\sum_{j\neq i} \big(1-z_j^2\big)^{1/2}\cos(\theta_i-\theta_j)  \right), & k = i,\\
2z_iz_k\big(1-z_i^2\big)^{-1/2}\big(1-z_{k}^2\big)^{-1/2}\cos(\theta_i-\theta_{k})+2, & k\neq i ,
\end{cases} \\
\frac{\partial^2 \hat H}{\partial \theta_k \partial z_i}   =   \begin{cases}
\displaystyle  2z_i\big(1-z_i^2\big)^{-1/2}\left(\sum_{j} \big(1-z_j^2\big)^{1/2}\sin(\theta_i-\theta_j)\right),& k = i,\\
 -2z_i\big(1-z_i^2\big)^{-1/2} \big(1-z_{k}^2\big)^{1/2}\sin(\theta_i-\theta_k),& k\neq i,
\end{cases} \\
\frac{\partial^2 \hat H}{\partial \theta_k \partial \theta_i}  =   \begin{cases}
\displaystyle -2\big(1-z_i^2\big)^{1/2}\left(\sum_{j\neq i} \big(1-z_j^2\big)^{1/2}\cos(\theta_i-\theta_j)\right), & k = i,\\
2\big(1-z_i^2\big)^{1/2}\big(1-z_k^2\big)^{1/2}\cos(\theta_i-\theta_k),& k\neq i,
\end{cases}
\\
\frac{\partial^2 \hat J}{\partial z_k \partial z_i}  =  \begin{cases}
-\big(1-z_i^2\big)^{-3/2}\big( z_{i-1}\big(1-z_{i+1}^2\big)^{1/2}\sin(\theta_{i+1}-\theta_i) &\\
\quad {} +z_{i+1}\big(1-z_{i-1}^2\big)^{1/2}\sin(\theta_i-\theta_{i-1})\big),   &\hspace*{-15mm} k = i, \\
\displaystyle\frac{-z_i}{\big(1-z_i^2\big)^{1/2}} \left( \big(1-z_{i+1}^2\big)^{1/2}\sin(\theta_{i+1}-\theta_i)- \displaystyle\frac{z_{i+1}z_{i-1}}{\big(1-z_{i-1}^2\big)^{1/2}}\sin(\theta_i-\theta_{i-1})   \right) &\\
 \quad {}  -\displaystyle\frac{z_{i-1}\big(1-z_{i+1}^2\big)^{1/2}}{\big(1-z_{i-1}^2\big)^{1/2}}\sin(\theta_{i-1}-\theta_{i+1}) ,                                                                                &\hspace*{-15mm} k = i-1,\\
\displaystyle\frac{-z_i}{(1-z_i^2)^{1/2}}\left( \frac{-z_{i-1}z_{i+1}}{\big(1-z_{i+1}^2\big)^{1/2}}\sin(\theta_{i+1}-\theta_i) + \big(1-z_{i-1}^2\big)^{1/2}\sin(\theta_i - \theta_{i-1})\right) & \\
\quad {}-  \displaystyle\frac{z_{i+1}\big(1-z_{i-1}^2\big)^{1/2}}{\big(1-z_{i+1}^2\big)^{1/2}}\sin(\theta_{i-1} - \theta_{i+1}),
&\hspace*{-15mm} k = i+1, \\
\end{cases}
\\
\frac{\partial^2 \hat J}{\partial \theta_k \partial z_i}   =  \begin{cases}
 -z_i\big(1-z_i^2\big)^{-1/2}\big( z_{i+1}\big(1-z_{i-1}^2\big)^{1/2}\cos(\theta_i-\theta_{i-1}) & \\
 \quad {} -z_{i-1}\big(1-z_{i+1}^2\big)^{1/2}\cos(\theta_{i+1}-\theta_i) \big), & k = i, \\
z_iz_{i+1}\big(1-z_{i-1}^2\big)^{1/2}\big(1-z_i^2\big)^{-1/2}\cos(\theta_i-\theta_{i-1}) & \\
\quad {}+ \big(1-z_{i+1}^2\big)^{1/2}\big(1-z_{i-1}^2\big)^{1/2}\cos(\theta_{i-1}-\theta_{i+1}) ,
  & k = i-1,\\
-z_iz_{i-1}\big(1-z_{i+1}^2\big)^{1/2}\big(1-z_i^2\big)^{-1/2}\cos(\theta_{i+1}-\theta_{i}) & \\
\quad {}- \big(1-z_{i+1}^2\big)^{1/2}\big(1-z_{i-1}^2\big)^{1/2}\cos(\theta_{i-1}-\theta_{i+1}),
  & k = i+1,
\end{cases}
\\
\frac{\partial^2 \hat J}{\partial \theta_k \partial \theta_i}   =  \begin{cases}
(1-z_i)^{1/2}\big( {-}z_{i+1}\big(1-z_{i-1}^2\big)^{1/2}\sin(\theta_i-\theta_{i-1})& \\
\quad{}-z_{i-1}\big(1-z_{i+1}^2\big)^{1/2}\sin(\theta_{i+1}-\theta_{i})  \big), &k = i, \\
z_{i+1}(1-z_i)^{1/2}\big(1-z_{i-1}^2\big)^{1/2}\sin(\theta_i-\theta_{i-1}),  & k = i-1,\\
 z_{i-1}(1-z_i)^{1/2}\big(1-z_{i+1}^2\big)^{1/2}\sin(\theta_{i+1}-\theta_{i}),  & k = i+1.
\end{cases}
\end{gather*}
Let
$ p = (0,s,0,s,\pi,-s)$ in cylindrical coordinates. The linearization of $X_{\hat J}$ at $p$ is
\begin{gather*}
A_{\hat J}(p) =\left(\begin{matrix}
0 & 1 & 1 & 0 & 0 & 0 \\
-1 & 0 & -1 & 0 & 0 & 0 \\
1 & -1 & 0 & 0 & 0 & 0 \\
0 & 0 & 0 & 0 & -1 & -1 \\
0 & 0 & 0 & 1 & 0 & 1 \\
0 & 0 & 0 & -1 & 1 & 0 \\
\end{matrix}\right).
\end{gather*}
The linearization of $X_{\hat H}$ at $p$ is
\begin{gather*}
A_{\hat H}(p) =2\left(\begin{matrix}
0 & 0 & 0 & 0 & -b^2 & b^2 \\
0 & 0 & 0 & -b^2 & 0 & b^2 \\
0 & 0 & 0 & b^2 & b^2 & -2b^2 \\
0 & b^{-2} & b^{-2} & 0 & 0 & 0 \\
b^{-2} & 0 & b^{-2} & 0 & 0 & 0 \\
b^{-2} &b^{-2} & 2b^{-2} & 0 & 0 & 0
\end{matrix}\right),
\end{gather*}
where $b^2 = 1-s^2$. These operators are independent and a quick computation shows that for any $-1< s < 1$ the operator on $L^{\perp}/L$ induced by $A_{\hat J}+A_{\hat H}$ has four distinct complex eigenvalues of the form~$A\pm iB$, $-A\pm iB$. Hence the critical point is rank 1 non-degenerate focus-focus.

\subsection*{Acknowledgements}

The author would like to thank his advisor Yael Karshon for her guidance and support and Leonid Polterovich for suggesting the study of integrable spin chains. Special thanks also goes to Joel Oakley for discussing his recent work on displaceability, Anton Izosimov for explaining his results on non-degenerate singularities, Peter Crooks for his editorial assistance, and the referees for their detailed feedback. The author was supported by NSERC PGS-D and OGS scholarships during the preparation of this work.

\pdfbookmark[1]{References}{ref}
\LastPageEnding

\end{document}